\documentclass{ifacconf}

\usepackage{natbib,graphics, graphicx, amsmath, float, amsfonts, amssymb, fp, tikz, pgfplots, mathrsfs, lipsum, mdframed, accents, booktabs, tabulary}

\newtheorem{proofpart}{Part}

\renewcommand{\vec}[1]{\mathbf{#1}}
\newcommand{\bi}{\begin{itemize}}
\newcommand{\ei}{\end{itemize}}
\newcommand{\beq}{\begin{enumerate}}
\newcommand{\eeq}{\end{enumerate}}
\renewcommand{\vec}[1]{\mathbf{#1}}
\newcommand{\be}{\begin{equation}}
\newcommand{\ee}{\end{equation}}

\newcommand{\bc}{\begin{center}}
\newcommand{\ec}{\end{center}}
\newcommand{\bd}{\begin{defn}}
\newcommand{\ed}{\end{defn}}
\newcommand{\bt}{\begin{thm}}
\newcommand{\et}{\end{thm}}
\newcommand{\bp}{\begin{proof}}
\newcommand{\ep}{\end{proof}}
\renewcommand{\l}{\left}
\renewcommand{\r}{\right}

\newcommand{\uhat}{\underaccent{\check}}

\DeclareMathOperator*{\argmax}{arg\,max}
\DeclareMathOperator*{\argmin}{arg\,min}

\begin{document}
\begin{frontmatter}

\title{Characterisation of Optimal Responses to Pulse Inputs in the Bergman Minimal Model}

\author[First]{Christopher Townsend}
\author[First]{Maria M. Seron}
\author[First]{Graham C. Goodwin}

\address[First]{School of Electrical Engineering and Computer Science, University of Newcastle, Australia (emails: chris.townsend@newcastle.edu.au, maria.seron@newcastle.edu.au, graham.goodwin@newcastle.edu.au)}

\begin{abstract}
	 The Bergman minimal model is a dynamic model of plasma glucose concentration. It has two input variables -- insulin delivery and carbohydrate intake. We investigate the behaviour of plasma glucose concentration predicted by the model given carbohydrate (CHO) inputs and commensurate insulin inputs. We observe that to maintain plasma glucose above a specified minimum concentration results in an unavoidable peak in plasma glucose. Additionally, we specify the timing and magnitude of a bolus pulse to minimise this unavoidable peak in plasma glucose concentration whilst attaining but not going below the desired minimum glucose concentration. Finally, we obtain necessary and sufficient conditions for the glucose concentration to be minimised.
\end{abstract}


\end{frontmatter}

\section{Introduction}

Diabetes is a chronic disease affecting over sixty million people \citep{wild04}. Diabetics are unable to correctly regulate blood glucose concentrations which, if not succesfully managed, leads to multiple adverse complications. Typically, management involves subcutaneous administration of insulin to minimise plasma glucose concentration whilst keeping it above a lower bound to avoid hypoglycaemia. Current treatment is invasive and often leads to poor outcomes. Hence, much recent effort has been devoted to developing an artificial pancreas which automates treatment \citep{harv10} and provides better control of glucose concentrations. The development of such systems and further treatment improvements requires an understanding of the dynamics of glucose regulation and pharmacokinetics of insulin. A number of models of glucose regulation have been proposed \citep{makr06}. One of these, the \emph{Bergman Minimal Model} (\cite{berg05, good15, kand09}), is a non-linear continuous-time model for glucose regulation. The model comprises a set of first order linear ordinary differential equations which govern the concentration and effectiveness of insulin:
\begin{align*}
  \frac{d}{dt} I_{sc}(t) &= -\frac{1}{\tau_1} I_{sc}(t) + \frac{1}{\tau_1}  \frac{ ID (t)}{C_l} \\
  \frac{d}{dt} I_p (t) &= -\frac{1}{\tau_2}  I_p (t) + \frac{1}{\tau_2} I_{sc} \\
  \frac{d}{dt} I_{eff}(t) &= - p_2 I_{eff}(t) +  p_2 S_I I_p (t)
\end{align*}%
\newline \quad \newline
and a non-linear ordinary differential equation which governs the plasma glucose concentration $g(t)$:
\begin{align*}
  \frac{dg}{dt}  = -g(t) \cdot  (I_{eff} (t) + G) + r(t) + E
\end{align*}%
where:
\begin{itemize}
 \item{$ID(t), I_{sc}(t), I_{p}(t)$ and $I_{eff}(t)$ -- are the delivery, subcutaneuos concentration, plasma concentration and insulin effectiveness, respectively. }
 \item{$\tau_1$ and $\tau_2$ -- are time constants.}
 \item{$C_l, S_I$ and $p_2$ -- are the clearance rate, insulin sensitivity and the insulin motility \citep{roy07}.}
 \item{$g(t)$ -- is the plasma glucose concentration.}
 \item{$E$ and $G$ -- are the endogenous glucose production and the effect of glucose on the uptake of plasma glucose and the suppression of endogenous glucose production, respectively.}
 \item{$r(t)$ -- is the glucose absorption from meals.}
\end{itemize}
A variety of physiological values for the above are derived from \citep{kand09} and given in Table 1 of \citep{good15}. For notational convenience, we rewrite the Bergman minimal model as the  following system of differential equations:
\begin{align}
\label{eq:eqs}
\begin{split}
  \dot z &= -d z + dk u \\
  \dot y &= -c y + cz \\
  \dot x &= -ax + aby\\
  \dot g &= - hg + w
 \end{split}
\end{align}%
where all variables and constants are positive, $u(t) = ID(t)$ is the input function,
\begin{align}
\label{eq:wandh}
\begin{split}
  h = x + G\\
  w = r + E
\end{split}
\end{align}%
   and the function $r$ is a given bounded function.

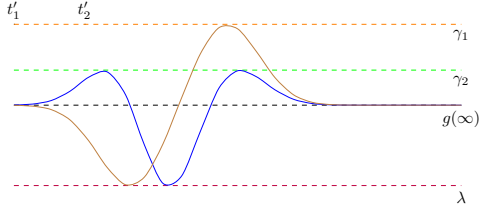
\begin{figure}[H]
	\begin{center}
\resizebox{0.35\textwidth}{!}{%
\begin{tikzpicture}
  
  \draw[.] (-2.5,3.73) -- (-2.5,3.73) node[above] {$t'_1$};
  \draw[.] (-1,3.73) -- (-1,3.73) node[above] {$t'_{2}$};

  \draw[-, dashed] (-2.5,2)  -- (7,2) node[below] {$g(\infty)$};
  \draw[-, dashed,orange] (-2.5,3.73) -- (7,3.73) node[below,black]  {$\gamma_1$} ;
  \draw[-, dashed,green] (-2.5,2.75)  -- (7,2.75) node[below,black] {$\gamma_2$};
  \draw[-, dashed,purple] (-2.5,0.28)  -- (7,0.28) node[below,black] {$\lambda$};

  \draw[scale=1,domain= -2.5:7,smooth,variable=\x, blue] plot ({\x},{2*(1.05*exp(-(\x)^2) - 1.6*exp(-(\x-0.61)^2))+ exp(-(\x-2)^2)+2});

 \draw[scale=1,domain= -2.5:7,smooth,variable=\x, brown] plot ({\x},{1.75*(exp(-(\x-2)^2) - exp(-(\x)^2))+ 2});
\end{tikzpicture}}
\caption{Glucose responses given a fixed function $w$ and two different injection times, where $t'_1$ corresponds to the injection time of the brown response and $t'_2$ the injection time of the blue response.}
\label{fig:opteg}
\end{center}
\end{figure}

   We develop necessary and sufficient conditions, given in Theorem \ref{thm:cream}, for the glucose response $g(t)$ to a pulse input function, $u(t)$, to be minimised. Additionally, these conditions give a non-linear version of the fundamental control limitation explored in Theorem 2 of \cite{medi15}, for the specific case of a single input pulse. In Figure \ref{fig:opteg} we show an example of the glucose responses to two pulse inputs delivered at times $t'_2 > t'_1$. The response which satisfies the conditions of Theorem \ref{thm:cream} has a lower maximum glucose concentration for the system whilst still maintaining the glucose concentration above a specified minimum concentration.

\section{Outline and Notation}
In Section \ref{sec:prelim} we state our assumptions and the constraints on the system. Additionally we prove some facts about the system and develop tools necessary for the subsequent Sections. In Section \ref{sec:opt} we prove necessary and sufficient conditions on the plasma glucose response to inputs $w$ and $u$ for the input $u$, comprising a single pulse, to be optimal. Finally, in Section \ref{sec:mult} we develop a sufficient condition for an input function, comprising multiple pulses, to be optimal.

\begin{table}[htbp]\caption{Notation}
\centering 
\begin{tabular}{r c p{5cm}}
\toprule
   
$\overline u, \hat u$ and $\hat U (\lambda)$ & -- & the basal input, the magnitude of the bolus input and an expression for $\hat u$ given in \eqref{eq:bolusbound}\\

$w,h$ and $g$ & -- & the bounded functions, see \eqref{eq:wandh}, and plasma glucose concentration, respectively.\\

    $\lambda$ and $\gamma$ & -- & the global minimum glucose concentration and the global maximum glucose concentration.\\

  $t', t_{\max}, t_{\min}$ and $\tau$ & -- & the delivery time, a time when the glucose concentration is at its global maximum, a time when the glucose concentration is at its minimum and the duration of the interval over which the bolus is delivered, respectively.\\

  $\mathcal{Y} (t)$ and $Y(t)$ & -- & the response of $x$ to the input $\hat u = 0$ and $\overline u =1$ and $\hat u =1$ and $\overline u =0$ respectively.\\

  $g(\infty)$ & -- & the steady-state glucose concentration $g(\infty) := \lim_{t \to \infty} g(t)$.\\

$g(h(u),w)$ & -- & the reponse of $g$ to the functions $h$ and $w$, where $h(u)$ is the response of $h$ to the input $u$.\\
\bottomrule
\end{tabular}
\label{tab:notation}
\end{table}

\section{Preliminaries, Assumptions and Constraints}
\label{sec:prelim}
   \subsection*{Assumptions}

Throughout we impose the following initial conditions: $z(0) = y(0) = k u(0)$, $x(0) = bk u(0)$ and $g(0) > 0$. We assume the function $w$ is positive and bounded. We also assume the input $u(t)$ is positive and bounded and  of the form:
\begin{equation}
	\label{eq:u}
	  u(t) = \bar{u}+ \hat{u} \chi_A (t)
\end{equation}
where the constant $\bar u$ is the \emph{basal} input, $\hat u$ is the magnitude of the \emph{bolus} input applied at some time $t'$, known as the \emph{delivery time}. The bolus input is held constant over the interval $A=[t',t' + \tau]$ and $\chi_A$ is the characteristic function of the interval $A$. The boundedness and positivity of $u(t)$ implies that $h$, given by \eqref{eq:eqs} and \eqref{eq:wandh}, is a continuous, positive and bounded function. We desire that there exist $\lambda > 0$ such that $g(t) \geq \lambda$ for all $t$. This is achieved if $\lambda$ is a global minimum of $g(t)$. We denote by $t_{min} \in \mathbb{R}_+$ a point such that $g(t_{min}) = \lambda$. Note, from \eqref{eq:eqs} that by setting $\dot g = 0$ at $t_{\min}$, we have $w(t_{\min} ) = \lambda h(t_{\min})$ at such $t_{\min}$.

\begin{defn}[Proper Input]
	For some $\lambda \leq g(0)$, an input function, $u(t)$, is \emph{proper}, if there exists $t_{\min}$ such that $g(h(u(t_{\min})),w) = \lambda$, $g(t) \geq \lambda$ for all $t$.
\end{defn}

The existence of a proper input, of the form \eqref{eq:u}, is established in Theorem \ref{thm:bolus}. Finally, unless otherwise stated we assume that $t_{\max} := \argmax_t {g(t)} < \infty$. The maximal time $t_{\max}$ exists as shown in Corollary \ref{cor:gammain}.

\subsection*{Bounds and System Properties}

\begin{lem}[Bounds]
\label{lem:bounds}
  Suppose $h$ and $w$ are bounded positive real-valued functionals, and $g$ is as in \eqref{eq:eqs}. Then there exist $\Gamma, \Lambda \in \mathbb{R}_+$ and constants $c_1$ and $c_2$ depending on the initial condition such that $\Gamma \geq \Lambda$ and:
  \[
     c_1 \exp\l (-\int h \r ) + \Lambda \leq  g(t) \leq c_2 \exp\l (-\int h \r ) + \Gamma
  \]
\end{lem}

\begin{pf}
A solution for $g$ is given by:
\[
  g(t) = \exp\l (-\int h \r ) \l ( c_3 + \int w \exp\l (\int h \r )\r)
\]%

where $c_3$ is the value of $g$ at the lower extreme of integration. Choose $\Gamma \in \mathbb{R}_+$ such that $w \leq \Gamma h$. Note such $\Gamma$ always exists since $w$ and $h$ are bounded positive functions. We obtain:
\begin{align*}
  g(t)& = \exp\l (-\int h \r ) \l ( c_3 + \int w \exp\l (\int h \r )\r) \\
      & \leq \exp\l (-\int h \r ) \l ( c_3 + \Gamma \int h \exp\l (\int h \r )\r) \\
      & = c_2 \exp\l (-\int h \r ) + \Gamma
\end{align*}

where $c_2 = c_3 -\Gamma$. Similarly, for $\Lambda$ such that $\Lambda h \leq w$, the lower bound on $g(t)$ is obtained and $c_1 = c_3 - \Lambda$.
\qed \end{pf}

\begin{rem}
\label{rem:upplow}
  The bounds, $\Gamma \geq \sup_t\l \{ \frac{w(t)}{h(t)} \r \}$ and $\Lambda \leq \inf_t \l \{ \frac{w(t)}{h(t)} \r \}$ in Lemma \ref{lem:bounds}, may be improved at any $t$ by taking a finite ordered partition $\mathcal{P} := \{t_0, \cdots, t_n, t\}$ of the interval $[0,t]$ and defining $\Gamma_i$ and $\Lambda_i$ such that $\Gamma_i h(t) \geq w(t)$ and $\Lambda_i h(t) \leq w(t)$ for all $t\in[t_i , t_{i+1}]$.
\end{rem}

\begin{lem}
\label{lem:tend}
Let $u(t)$ be as in \eqref{eq:u} and choose $\lambda \leq g(0)$. Then:

\begin{enumerate}
\item{ The response $x(t)$ is separable into its basal and bolus responses i.e.
  \begin{equation}
    \label{eq:x}
    x(t) = x(\bar u) + x(\hat u) := \bar u \mathcal{Y}(t) + \hat u Y(t)
  \end{equation}
  where $\mathcal{Y}(t) $ is the response of $x$ to the input with $\bar u = 1$, $\hat u =0$ and $Y(t)$ is the response of $x$ to the input~\eqref{eq:u} with $\bar u = 0$, $\hat u=1$.}
  
  \item{Under the assumed initial conditions, $x(\overline u) = bk\overline u$ for all $t$ i.e. $\mathcal{Y} (t) = bk$. Furthermore, if $A$ is bounded. Then:
  \[
    \lim_{t \to \infty} x(t) = bk \overline u 
  \]}
  
  \item{if $A$ is empty and the basal input $\overline u \leq \frac{E - \lambda G}{\lambda bk}$ then $g(t) \geq \lambda$ for all $t$.}
 \end{enumerate}
\end{lem}

\begin{pf}

\begin{proofpart}
  This follows by linearity of the $(z,y,x)$-system in~\eqref{eq:eqs}.
\end{proofpart}

\begin{proofpart}
Follows from the solutions to the first order linear ordinary differential equations in \eqref{eq:eqs} and $\mathcal{Y}(t)= bk$ from the specified initial conditions.
\end{proofpart}

\begin{proofpart}
By Lemma \ref{lem:bounds} there exists $\Lambda \in (0, g(0)]$ satisfying $\Lambda h \leq w$, which guarantees that $g(t) \geq \Lambda$ for all $t$. As, from \eqref{eq:wandh}, $w(t) \geq E$ and $h(t) = h(\overline u) = x(\overline u) + G \leq \overline u bk +G$ by Lemma \ref{lem:tend}, such $\Lambda$ may be chosen to satisfy the inequality:
  \[
    \Lambda \leq \frac{E}{bk \overline u + G}
  \]%
  Therefore, to ensure that $\Lambda \geq \lambda$ we require:
  \begin{align*}
    \lambda \leq \frac{E}{bk \overline u + G}
    \iff \overline{u} \leq \frac{E- \lambda G}{\lambda bk}
  \end{align*} 
  \end{proofpart}
  \hspace{8cm} \qed
  \end{pf}

\begin{rem}
\label{rem:hu}
  We define the \emph{steady-state} of $g$ to be $g(\infty):=\lim_{t \to \infty} g(t)$, when $\lim_{t \to \infty}Y(t)  = 0$ and $\lim_{t \to \infty} w(t) = E$ i.e. it is the limit of the response of $g(t)$ when the only input is the constant input $\overline u$. The steady-state may be set to be any positive real number. Since, if $\overline u$ is constant then $\dot g \to 0$ as $t \to \infty$. Indeed, assuming that $A$ is bounded, setting:
  \begin{align}
  \label{eq:basal}
    \overline u = \frac{1}{kb}\l (\frac{E}{g(\infty)} - G \r)
  \end{align}%
  gives the result, where $g(\infty) \geq \lambda$ is some specified value. In this case: 
  \[
    h(\overline u) :=  x(\overline u) + G \leq \frac{E}{g(\infty)}
  \]%
  We note, the constant, $c_3$, from Lemma \ref{lem:bounds} corresponds to the initial plasma glucose concentration, $g(0)$. Henceforth, we fix $g(\infty) := g(0)$ and $\overline u$ to be as in~\eqref{eq:basal}.
\end{rem}

The nature of the system dynamics \eqref{eq:eqs}--\eqref{eq:wandh} and the positivity of the inputs induces a monotonic relationship between the insulin input and glucose concentration. This property is proven in the intermediate Lemma \ref{lem:forprooftowork}. Note, we define $u_1 (t) > u_2  (t)$ if there exists $s$ such that $u_1(s) > u_2(s)$ and $u_1 (t) \geq u_2(t)$ for all $t$.

\begin{lem}
	\label{lem:forprooftowork}
	Suppose $w$ is fixed. Then $g(t)$ is a strictly monotone function of the input $u(t)$.
\end{lem}

\begin{pf}
    Fix $w$ and let $u_1$ and $u_2$ be two inputs with delivery time $t'$ such that $\hat u_1 < \hat u_2$. Denote by $h_1, g_1$ and $h_2$ and $g_2$ their respective responses. Since $h$ is a monotone function of the input $u$ we have that $h_1 < h_2$ for all $t > t'$. 
A solution for $g(t)$ for $t \geq t'$ is given by:
\begin{align*}
g(t) & = g(t')\exp\l (-\int_{t'} ^t h(s) \, ds \r ) \\ & \quad + \int_{t'} ^t w(s) \exp\l (- \int_{s} ^t h(\xi) \, d \xi \r ) \, ds
\end{align*}%
As $g_1 (t') = g_2 (t')$, because the inputs are identical before $t'$, and:
\begin{align*}
  \exp\l (-\int_{t'} ^t h_1 (s) \, ds \r ) > \exp\l (-\int_{t'} ^t h_2 (s) \, ds \r ) \\
  \exp\l (- \int_{l} ^t h_1(\xi) \, d \xi \r ) \geq \exp\l (- \int_{l} ^t h_2(\xi) \, d \xi \r )
\end{align*}%
for all $t > t'$ and $l \leq t$, we have that $g_1 (t) > g_2(t)$ for all $t> t'$.
\qed \end{pf}

Theorem \ref{thm:bolus} proves the existence of a bolus input delivered at any $t'$ which achieves a specified minimum $\lambda > 0$ and thus proves the existence of proper inputs of the form \eqref{eq:u}.

\begin{thm}[Insulin Bolus]
\label{thm:bolus}
  Suppose $u(t)$ is of the form \eqref{eq:u}. Fix $\tau$ and  $t'$ -- the input time i.e. $A:= [t', t'+ \tau]$, choose $\lambda \in (0, g(t')]$ and suppose $\bar u$ is as in Remark \ref{rem:hu}. Then there exists $\hat u$ such that $u(t)$ is proper.
\end{thm}

\begin{pf}
	Denote by $g(\hat u)$ the response of $g(t)$ to the input $u(t) := \overline u + \hat u \chi_{A}$. By Lemmas \ref{lem:bounds} and \ref{lem:forprooftowork}, there exist $\hat o$ and $\hat v$ such that $\min_{t\geq t'} g(\hat o) \geq \lambda$ and $\min_{t\geq t'} g(\hat v) \leq \lambda$. 

Suppose $\min_{t\geq t'} g(\hat v) < \min_{t\geq t'} g(\hat o)$. We recursively define the sequences $\overline o := (\hat o _i)_{i=0} ^\infty$ and $\overline v := (\hat v _i)_{i=0} ^\infty$ by $\hat o _0 = \hat o$ and $\hat v_0 = \hat v$ and $\hat o_i$ the greatest element of the following finite ordered partition of the interval $[\hat o_{i-1}, \hat v_{i-1}]$:

\begin{align*}
	{L}_i  := &\l \{\hat o_{i-1}, \frac{ (n-1) \hat o_{i-1} + \hat v_{i-1}}{n}, \cdots, \r . \\  &\quad \cdots,  \l. \frac{k_i \hat o_{i-1} + (n-k_i) \hat v_{i-1}}{n}, \cdots, \hat v_{i-1} \r \}
 \end{align*}%
 where $n \in \mathbb{N}$ is arbitary and $k_i \leq n$, such that the response:
 \[ 
  g \l (\frac{k_i \hat o_{i-1} + (n-k_i) \hat v_{i-1}}{n} \r ) \geq \lambda
  \]%
 for all $t \geq t'$. Similarly, $\hat v_i$ is defined to be the least element of $L_i$ such that, for all $t \geq t'$:
 \[ 
    g \l (\frac{k_j \hat o_{i-1} + (n-k_j) \hat v_{i-1}}{n} \r ) \leq \lambda
 \]%
 The sequence $\overline o$ is a monotone increasing sequence bounded above by $\hat v_i$ for all $\hat v_i \in \overline v$ and therefore has a limit $o$. Similarly, $\overline v$ is a monotone decreasing sequence bounded below by $\hat o_i$ for all $\hat o_i \in \overline o$ and thus has a limit $v$. It remains to show that these two limits are equal. If either sequence is eventually constant then both are constant and equal. As either sequence is constant only if $\min_{t \geq t'} g(t) = \lambda$, in which case, by construction of the sequences, both sequences would have the same value. Suppose, instead, for all $i$ that $\hat o_i < \hat v_i$. We see that if:
 \[ 
  \hat o _{i+1} = \frac{k_i \hat o_{i} + (n-k_{i})\hat v_i}{n}
 \]%
 Then, by Lemma \ref{lem:forprooftowork}, $\hat v_{i+1}$ must be the next element of $L_i$, that is:
 \[ 
  \hat v_{i+1} = \frac{(k_i-1) \hat o_{i} + (n-k_{i}+1)\hat v_i}{n}
 \]%
 Thus:
 \begin{align*}
  \hat v_{i+1} - \hat o_{i+1} &= \frac{1}{n} \l (\hat v_i - \hat o_i \r) \\
                              &\qquad \vdots\\
                              &= \frac{1}{n^{i+1}} \l (\hat v_0 - \hat o_0 \r)
 \end{align*}%
 
 i.e. $\lim_{i \to \infty} \l (\hat v_{i+1} - \hat o_{i+1} \r ) = 0$ i.e. $v = o$. Thus $\lambda \leq \min_{t' \geq t} g(v) = \min_{t' \geq t} g(o) \geq \lambda$. Setting $\hat u = v$ gives $\min_{t' \geq t} g(\hat u) = \lambda$ and $g(\hat u) \geq \lambda$ for all $t\geq t'$.
\qed \end{pf}

Corollary \ref{cor:bolus} provides an explicit expression for the magnitude of the bolus input which achieves the gloal minimum of the glucose concentration.

\begin{cor}[Insulin Bolus Bound]
\label{cor:bolus}
Fix $t'$ and choose $\lambda > 0$ and let $u= \overline{u} + \hat{u} \chi_A$, see \eqref{eq:u}, where $\bar u$ is as in Remark \ref{rem:hu}. Suppose the input $\hat u$ is as in Theorem \ref{thm:bolus}. Then the input satisfies, $\hat u \leq \hat U (\lambda)$, where:
  \begin{align}
  \label{eq:bolusbound}
    \hat U (\lambda) := \l ( \frac{w(t_{min})}{\lambda} - G - x(\overline u, t_{\min})\r) \l( \frac{1}{Y(t_{min})} \r)
  \end{align}
  and $g(t) \geq \lambda$ for all $t$. In particular, if $\hat u = \hat U(\lambda)$ then $g(t_{\min}) = \lambda$. 
\end{cor}

\begin{pf}
  We have the following: 
  \begin{align}
	  \begin{split}
    \label{eq:2}
    \dot g(t_{min}) = & -g(t_{min}) h(t_{min})+w(t_{min}) = 0 \\ &  \iff g(t_{min}) = \frac{w(t_{min})}{h(t_{min})}
    \end{split}
  \end{align}%
  Suppose $g(t_{min}) < \lambda$. Then, from \eqref{eq:wandh}, \eqref{eq:x}, \eqref{eq:2}, and as $h(t) > 0$:  
  \begin{align*}
	  w(t_{min}) & < \lambda h(t_{min})= \lambda (x(\bar u, t_{\min}) + \hat u Y(t_{min}) +G) \\ & \iff \hat{u} > \hat{U}(\lambda).
  \end{align*}
Thus, \eqref{eq:bolusbound} implies $g(t) \ge \lambda$ for all $t$ and the result follows.
\qed \end{pf}

Corollary \ref{cor:gammain} proves the existence of a finite global maximum for glucose concentration responses to proper inputs.

\begin{cor}[An Upper Bound]
\label{cor:gammain}
  Choose $\lambda \leq g(0)$. Suppose $u(t)$ is proper. Then there exists $t_{\max} \in \mathbb{R}_+ ^* := \mathbb{R}_+ \cup \{\infty \}$ such that $g(t) \leq g(t_{\max}) =: \gamma$ for all $t$ and $\gamma = \lambda$ if and only if $g(t) = \lambda$ for all $t$. Furthermore, if $t_{\max} < \infty$. Then:
    \begin{align}
      \label{eq:maxgam}
      \gamma = \frac{w(t_{\max})}{h({t_{\max}})} = \frac{\alpha_1 \lambda}{\alpha_2 + \alpha_3 \lambda}
    \end{align} 
    where $\alpha_1 := w(t_{\max})$, $\alpha_2 :=  w(t_{\min}) \l (\frac{Y(t_{\max})}{Y(t_{\min})} \r)$,  $\alpha_3 := \l (G+ x(\overline u, t_{\max}) \r) \l( 1 -  \frac{Y(t_{\max})}{Y(t_{\min})} + x(\overline u, t_{\max}) - x(\overline u, t_{\min})\r)$.
\end{cor}

\begin{pf}

    If there is $s \in \{t : \dot g(t) = 0\}$ such that $g(s) \geq g(t)$ for all $t$. Then $t_{\max} = s$. Otherwise $g(t)$ must increase as $t \to \infty$. We may take a monotone increasing sequence $(g(t_i))_{i=0} ^\infty$, where $t_0 \geq 0$ and $g(t_k) \in (g(t_i))_{i=0} ^\infty$ only if $g(t_k) > g(t)$ for all $t \in [0, t_k)$ i.e. $g(t_k)$ is the \emph{peak} of the function $g$. By Lemma \ref{lem:bounds} the sequence, $(g(t_i))_{i=0} ^\infty$, is bounded above and thus converges to $\overline g \leq \Gamma$. By construction $g(t) \leq \overline g$ for all $t$. We see that $\gamma = \lambda$ if and only if $g(t) = \lambda$ for all $t$ follows by definition of $\gamma$.
  
  Suppose $t_{\max} < \infty$. Then \eqref{eq:maxgam} follows from rearranging the differential equation $\dot g = -gh + w$ evaluated at $t_{\max}$ and substituing in the formula for $\hat u$ given by \eqref{eq:bolusbound}.
\qed \end{pf}

Corollary \ref{cor:seesaw} shows that, the higher the minimum glucose concentration, the higher the peak glucose concentration.

\begin{cor}
\label{cor:seesaw}
Choose $\lambda < \lambda' \leq g(0)$. Let $u(t)$ and $u'(t)$ be inputs, of the form \eqref{eq:u}, with common delivery time $t'$, which are proper for $\lambda$ and $\lambda'$ respectively. Then $\gamma < \gamma'$.
\end{cor}

\begin{pf}
	Note $\gamma' \geq \lambda'$ and $\gamma > \lambda$. Denote by $t_{min}$ and $t'_{min}$ the times at which $g(t) = \lambda$ and $g'(t) =\lambda'$ respectively. If for example $h' > h$ for some $t>t'$. Then $h' > h$ for all $t > t'$. This is because $u(t)$ and $u'(t)$ are of the form \eqref{eq:u} and have common delivery time. Suppose that $h' \geq h$ for all $t \geq t'$. This implies that $g' \leq g$ for all $t \geq t'$. In particular at $t_{\min} \ge t'$ we have that $\lambda = g(t_{\min}) \geq g'(t_{\min}) \geq \lambda'$ contradicting $\lambda < \lambda'$. Thus $h > h'$ for all $t > t'$. Finally, as $h$ is a monontone function of the input $u(t)$, we see that $\hat u > \hat u'$. This implies that $g(t) < g'(t)$ for all $ t > t'$.
\qed \end{pf}

\section{Optimal Inputs}
\label{sec:opt}

In this section we give necessary and sufficient conditions on the delivery time of a proper input such that the glucose response is optimal i.e. the maximum glucose concentration $\gamma$ is minimised. Throughout, we fix the length, $\tau$, of the interval over which the bolus is delivered. This ensures that the function $Y(t - t')$ in \eqref{eq:x} is invariant under translation by $t'$ -- the delivery time. Under these conditions we establish a property of the function $h$, defined in \eqref{eq:wandh}, in  Lemma \ref{lem:morelemmas}.

\begin{lem}
\label{lem:morelemmas}
For any two distinct proper inputs $u_1$ and $u_2$, delivered at times $t'_1$ and $t'_2$ respectively, with responses $h_1$ and $h_2$, there exists $t_i$ such that either:
\begin{align}
\label{eq:firstcase}
   \l \{ \begin{array}{lc} 
        h_1 > h_2,&  t \in (\min \{t'_1, t'_2\}, t_i) \\ 
        h_1 = h_2,&  t = t_i \\ 
        h_1 < h_2,&  t > t_i 
  \end{array}\r .
\end{align}%
Or:
\begin{align}
\label{eq:secondcase}
  \l \{ \begin{array}{lc} 
        h_1 < h_2,&  t \in (\min \{t'_1, t'_2\}, t_i) \\ 
        h_1 = h_2,&  t = t_i \\ 
        h_1 > h_2,&  t > t_i 
  \end{array}\r .
\end{align}%
\end{lem}

\begin{pf}
	Indeed, if $h_2 (t) > h_1(t)$ or $h_2 (t) < h_1(t)$ for all $t > \min \{t'_1, t'_2\}$ then, by Lemma \ref{lem:forprooftowork} there would exist $t$ such that $g_2 (t) <  \lambda$ or $g_1 (t) < \lambda$. Implying that either $\hat u _1$ or $\hat u_2$ are not proper. 
\qed \end{pf}

In the subsequent proofs of Lemmas \ref{lem:anotherone}, \ref{lem:optimaltime} and \ref{lem:optglu}, we only present the first case, \eqref{eq:firstcase}, as the other case, \eqref{eq:secondcase}, follows by a similar argument.

\begin{rem}
	\label{rem:labelit}
 For fixed $w$ the functions $h_1$ and $h_2$ satisfy \eqref{eq:firstcase} only if $t'_1 < t'_2$ i.e. only if $\hat u_1$ is delivered before $\hat u_2$ as the function $Y(t)$ is independent of the magnitude $\hat u$.
\end{rem}

We define the optimal, minimised, glucose response as the input control strategy that ensures that the maximum glucose concentration is minimised given the control and system limitations. Formally:

\begin{defn}[Minimised Response]
We say the response $g(t) = g(w,h)$ is \emph{minimised} if  $\max_t g(w,h') > \max_t g(w,h)$, for all $h' \neq h$,
\end{defn}

\begin{defn}[Optimal Delivery Time]
  We say a delivery time $t'$ of a bolus input $\hat u \chi_{[t',t' + \tau]}$, where $\hat u$ is given by (\ref{eq:bolusbound}), is \emph{optimal} if the response $g(t)$ is minimised.
\end{defn}

\begin{lem}
\label{lem:anothernotherone}
Suppose $w(t)$ is a continuous and bounded positive functional and $h(t)$ is the response to a proper input of the form \eqref{eq:u}. Then, for any $\varepsilon \in [0, g(0) - \lambda)$, there exist at most finitely many $t$ such that the response $g(t) < \lambda + \varepsilon$ and $\dot g(t) = 0$.
\end{lem}

\begin{pf}
  This follows as $w(t)$ is bounded below and the response $Y(t)$ is continuous and converges to $0$.
\qed \end{pf}

\begin{lem}
\label{lem:anotherone}
	Suppose $g (h_1, w) = g_1(t)$ is a response to a proper input with  bolus $\hat u_1$ delivered at time $t'_1$ such that there is a unique minimum $t_{1,\min}$ i.e. $g_1 (t_{1,\min}) = \lambda$ and $g_1 (t) > \lambda$ for all $t \neq t_{1,\min}$. Then there exists a proper bolus input $\hat u_2$ delivered at $t'_2$ and a time $t_i \geq \max\{t'_1,t'_2\}$ at which $h_1(t_i) = h_2(t_i)$ such that the response $g_2(h_2, w) := g_2(t)$ attains its minimum at $\lambda$ and satisfies:
\[
  \l \{ \begin{array}{lc} 
        g_1 < g_2,&  t < t_g \\ 
        g_1 = g_2,&  t = t_g \\ 
        g_1 > g_2,&  t > t_g 
  \end{array}\r .
\]%
for some time $t_g \in [t_i, t_{1, \max})$, where $t_{1,\max} := \min \{s > t_{1,\min} : g_1(s) \geq g_1(t) \, \forall \, $t$\}$.
\end{lem}

\begin{pf}
	As $g$ is a continuous function of $h$, for all $\varepsilon > 0$ we may find $\delta > 0$ such that $|h_1 - h_2| < \delta$ implies $|g_2 - g_1| < \varepsilon$, for all $t$. Such $h_2$ exists and is of the form $x(t) + G$, where $x(t)$ is as in \eqref{eq:x}, as $x(t)$ is a continuous function of the input $u(t)$. Thus for all $\delta > 0$ there exists $\delta' > 0$ such that $\max |u_1  - u_2 | < \delta'$ implies $|h_1 - h_2| < \delta$. Furthermore, by Theorem \ref{thm:bolus}, we may assume $u_2$ is proper and of the form \eqref{eq:u}. Thus, by Lemma \ref{lem:morelemmas} and Remark \ref{rem:labelit}, $u_2$ may be chosen such that there exists $t_i$ for which:
\[ \l \{ \begin{array}{lc} 
        h_1 > h_2,&  t \in (t_1 ', t_i) \\ 
        h_1 = h_2,&  t = t_i \\ 
        h_1 < h_2,&  t > t_i 
  \end{array}\r .
\]%
 Define $\overline \lambda := \min\{g_1(t) : t \neq t_{1,min} \wedge \dot g_1 (t) = 0\}$ or $g(0)$ if this minimum does not exist. Such $\overline \lambda > \lambda$ exists by Lemma \ref{lem:anothernotherone} and $\overline \lambda \leq g(0)$, by construction.
 Choosing $\varepsilon < \min\{\gamma_1 - g(0), \overline \lambda - \lambda\}$ or $\varepsilon < \overline \lambda - \lambda$ if $\gamma_1 - g(0)= 0$, where $\gamma_1 := \max(g_1(t))$, implies that $t_{2,min} < t_{1,\max}$. Note that $t_{2, \min} > t_i$ as if it were not there would exist $t$ such that $g_1 (t) < \lambda$ since $g_1 < g_2$ for all $t \in (t_1', t_i)$. Also, $g_2(t_{2, \min}) = \lambda < g_1(t_{2,\min})$, by assumption. By the Intermediate Value Theorem there is a $t_g \in [t_i, t_{1, \max} )$ such that $g_1 (t_g) = g_2 (t_g)$. As $h_2 > h_1$ for all $t > t_i$ we see that $g_2 (t) < g_1 (t)$ for all $t > t_g$.
\qed \end{pf}

\begin{rem}
  Similarly to Lemma \ref{lem:anotherone}, we may show that if $t_{1, \max} < t_{1, \min}$ then there exists $h_2$ with response $g_2$ and $t_g \in [t_i,t_{1, \min})$ such that:
\[
  \l \{ \begin{array}{lc} 
        g_1 > g_2,&  t < t_g \\ 
        g_1 = g_2,&  t = t_g \\ 
        g_1 < g_2,&  t > t_g 
  \end{array}\r .
\]%
\end{rem}

\begin{lem}[Single Minimum]
\label{lem:optimaltime}
  Suppose $g(h,w) = g(t)$ is a response, to a proper input $u$, for which there is a unique $t_{\min}$ such that $g(t_{\min}) = \lambda$. Then $g(t)$ is minimised if and only if $\max_{t < t_{\min}} g(t) = \max_{t > t_{\min}} g(t)$.
\end{lem}

\begin{pf}
Without loss of generality, suppose that for some $h_1$ we have that:
\[
  \overline \gamma_1 := \max_{t < t_{1,\min}} g_1(t) < \max_{t > t_{1,\min}} g_1 (t) =: \gamma_1
\]%
where $g_1 (t) := g(h_1,w)$. The existence of a unique $t_{1,\min}$ implies there is a proper input with non-zero bolus $\hat u_1$ delivered at some time $t'_1$. Define $\bar t _1 := \argmax_{t < t_{1,\min}} g_1(t)$ and $\hat t_1 := \argmax_{t > t_{1,\min}} g_1 (t)$. By Lemma \ref{lem:anotherone} there exists $h_2$ such that:
\[
  \l \{ \begin{array}{lc} 
        g_1 < g_2,&  t < t_g \\ 
        g_1 = g_2,&  t = t_g \\ 
        g_1 > g_2,&  t > t_g 
  \end{array}\r .
\]%
Proceeding as in the proof of Lemma \ref{lem:anotherone}, choosing $\varepsilon < \min\{\gamma_1 - g(\infty), \overline \lambda - \lambda, \gamma_1 - \overline \gamma_1\}$ implies that $g_2 < \gamma_1$ for all $t$, since the choice of $\varepsilon$ ensures that $t_g < \hat t_{1}$.

Now, suppose that $\max_{t < t_{1,\min}} g_1(t) = \max_{t > t_{1,\min}} g_1 (t) := \gamma_1$ and $g_1$ is not minimised. We observe that for $g_1$ not to be minimised there must exist $h_2$ with response $g_2$ such that $g_2 (\bar t _1) < \gamma_1$ and $g_2 (\hat t_1 ) < \gamma_1$, which implies $t_g < \overline t_1$. So $g_2 (t_{1,\min}) < g_1 (t_{1,\min}) = \lambda$. This condradicts the constraint on $g_2$. Hence no such $h_2$ exists.
\qed \end{pf}

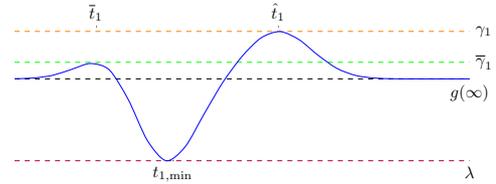
\begin{figure}[H]
	\begin{center}
\resizebox{0.35\textwidth}{!}{%
\begin{tikzpicture}
  
  \draw[.] (-0.8,3.1) -- (-0.8,3.1) node[above] {$\overline{t}_1$};
  \draw[.] (3,3.1) -- (3,3.1) node[above] {$\hat t_{1}$};
  \draw[.] (0.8,0.28) -- (0.8,0.28) node[below] {$t_{1,\min}$};

  \draw[-, dashed] (-2.5,2)  -- (7,2) node[below] {$g(\infty)$};
  \draw[-, dashed,orange] (-2.5,3.0) -- (7,3.0) node[right,black]  {$\gamma_1$} ;
  \draw[-, dashed,green] (-2.5,2.35)  -- (7,2.35) node[right,black] {$\overline \gamma_1$};
  \draw[-, dashed,purple] (-2.5,0.28)  -- (7,0.28) node[below,black] {$\lambda$};

  \draw[scale=1,domain= -2.5:7,smooth,variable=\x, blue] plot ({\x},{2*(1.05*exp(-(\x)^2) - 1.65*exp(-(\x-0.4)^2))+ exp(-(\x-3)^2)+2});

\end{tikzpicture}}
\caption{Glucose response for functions $w$ and $h_1$ showing two unequal maxima about a single minimum, where $\overline \gamma_1 := \max_{t < t_{1,\min}} g_1(t)$ and $\gamma_1 := \max_{t > t_{1,\min}} g_1(t)$.}
\label{fig:unequalglu}
\end{center}
\end{figure}

\begin{lem}[Single Maximum]
\label{lem:optglu}
  Suppose $g(h,w) = g(t)$ is a response, to a proper input, for which there exist distinct $\underline{t}$ and $\uhat t$ such that $g(\underline{t}) = g(\uhat t) = \lambda$ but a single $t_{\max} := \argmax\{g(t)\}$. Then $g(t)$ is minimised if and only if $\argmax_t \{g(t)\} \in (\underline{t},\uhat t)$.
\end{lem}

\vspace{10pt}

\begin{pf}
  Suppose, for some $h_1$, that:
  \[ 
    \min_{t < t_{1, \max}} \{ g_1 (t) \} = \min_{t > t_{1, \max}} \{ g_1(t) \} = \lambda
    \] Define $\underline{t}_1 := \min \{t < t_{1, \max} :  g_1 (t) = \lambda \}$ and ${\uhat{t}}_1  := \max \{t > t_{1, \max} :  g_1(t) = \lambda \}$. Suppose $h_2 \neq h_1$ is a response to a proper input $u_2$ as in Lemma \ref{lem:anotherone}. As $g_2 \geq \lambda$ for all $t$, the crossing time $t_g$ satisfies either $t_g < \underline t_1$ or $t_g > \uhat{t}_1$. In both cases $g_2(t) > g_1(t)$ for all $t \in [\underline t_1, \uhat{t}_1]$ which implies $\max\{g_2(t)\} > \max\{g_1 (t) \}$.

  Now, suppose that there are at least two distinct $t$ such that $g_1(t) = \lambda$ and, without loss of generality, that:
   \[
     \min_{t < t_{1, \max}} \{ g_1(t)\} < \min_{t > t_{1, \max}} \{g_1(t)\}
   \]%

   As in Lemma \ref{lem:anotherone} there exists $h_2$ and $t_i$ with response $g_2$ such that:
   \[
  \l \{ \begin{array}{lc} 
        g_1 < g_2,&  t < t_g \\ 
        g_1 = g_2,&  t = t_g \\ 
        g_1 > g_2,&  t > t_g 
  \end{array}\r .
\]%
where $t_g \in [t_i, t_{1, \max})$. In this case, proceeding as in Lemma \ref{lem:optimaltime}, there exists $h_2$ such that $g_2 (t) < \gamma_1$, where $\gamma_1 := \max\{g_1(t)\}$. 
  \qed  \end{pf}

Lemmas \ref{lem:optimaltime} and \ref{lem:optglu} show that for a given $w$ and input of the form $u(t) := \overline u + \hat u \chi_{[t', t'+ \tau]}$, see \eqref{eq:u}, the maximum glucose concentration $g(t_{\max})$ is minimised if and only if the maximum occurs between two minima, where $g(t) = \lambda$, or the minimum between to equal maxima. We state this formally in Theorem \ref{thm:cream}.

\begin{thm}[Multiple Extrema]
\label{thm:cream}
  Suppose $g(h,w) = g(t)$ is a response to an input of the form \eqref{eq:u}. Then $g(t)$ is minimised if and only if there exists $t_{\min}$ such that $g(t_{\min}) = \lambda$ and either: $\argmax_t \{g(t)\} \in (\underline{t},\uhat t)$ where $\underline t, \uhat t \in \argmin\{g(t)\}$ or there is $\uhat t \in \argmin\{g(t)\}$ such that $\max_{t< \uhat t} \{g(t) \}= \max_{t< \uhat t} \{g(t) \}$.
\end{thm}

	 \section{Multiple Pulses}
	 \label{sec:mult}
 The optimality conditions given in Theorem \ref{thm:cream} apply to an input $u(t)$ with only a single bolus input delivered at some $t'\in \mathbb{R}_+$. We provide a sufficient condition for optimality of an input with finitely many bolus inputs i.e. we consider inputs of the form:
\begin{align}
	\label{eq:genu}
		u(t,T) := \overline u + \sum_{i=0} ^N \hat u_i \chi_{[t'_i, t'_i + \tau]}
\end{align}%
where $T := (t_i)_{i=0} ^N$ is a finite sequence of delivery times and the magnitude $\hat u_i$ of each bolus input is given sequentially by Theorem \ref{thm:bolus}\footnote{Theorem \ref{thm:bolus} does not explicitly need the input to have the form~\eqref{eq:u}}, if there exists $t > t'_i$ such that $g(t) > g(0)$ or if $i = 0$. Otherwise $\hat u_i := 0$. This ensures that each bolus input, and therefore the input, $u(t,T)$, is proper and avoids unnecessary inputs. Thus we may assume that each $\hat u_i > 0 $. By abuse of notation, we denote $u(t,N) = u(t,T)$, where $N$ is the length of the sequence $T$.

\begin{lem}
	\label{lem:extend}
	Suppose $u_1(t,N)$ and $u_2(t,N)$ are distinct inputs of the form \eqref{eq:genu}. Then there exist at most $2N$, $t_{g,i}$ such that $g_1 (t_{g,i}) = g_2(t_{g,i})$ and for which one of the following, with either direction or order of the inequalities, is satisfied:
 \[
  \l \{ \begin{array}{lc} 
		  g_1 < g_2,&  t \in (t_{g,i-1}, t_{g,i}) \\ 
		  g_1 = g_2,&  t = t_{g,i} \\ 
		  g_1 \geq g_2,&  t \in (t_{g,i}, t_{g,i+1}) 
  \end{array}\r .
\]%
Or:
\[
  \l \{ \begin{array}{lc} 
		  g_1 < g_2,&  t \in (t_{g,i-1}, t_{g,i}) \\ 
		  g_1 = g_2,&  t = t_{g,i} \\ 
		  g_1 < g_2,&  t \in (t_{g,i}, t_{g,i+1}) 
  \end{array}\r .
\]%
where $t_{g,N+1} := \infty$. Furthermore, $t_{g,0} := \inf\{t : g_1 (t) \neq g_2(t) \}$ must exist.
\end{lem}

\begin{pf}
	This follows as $g$ is a monotonic function of $h$ which is a monotonic function of the input $u$.
\qed \end{pf}

\begin{defn}
	Let $u_1(t,N)$ and $u_2(t,N)$ be two distinct inputs of the form \eqref{eq:genu}. The points $\{t_{g,0}, t_{g,1}, \cdots, t_{g,2N -1} \}$ defined in Lemma \ref{lem:extend} are the \emph{intersection points} of the responses $g(h(u_1(t,N)),w) = g_1(t)$ and  $g(h(u_2(t,N)),w) = g_2(t)$.
\end{defn}

	 \begin{thm}
		 \label{thm:mult}
		 Suppose $u(t,N)$ is proper, the sum of the number of global maxima and minima, of the response $g(h(u(t,N)),w)$, is $2N+1$ and these minima and maxima are interlaced. Then $u(t,N)$ is optimal, over inputs of the same $N$.
	 \end{thm}

	 \begin{pf}
		 Suppose $g_1$ is a response as in the statement of the Theorem and that there is an input $u_2(t,N) \neq u_1(t,N)$ such that $\max \{ g_2 (t) \} < \max \{ g_1(t)\}$. Additionally, assume that $t_{g,i} \neq t_{1,i,\min}$, where $t_{1,i,\min}$ is the $i^{\text{th}}$ minimum of $g_1$ i.e. the intersection points of $g_1$ and $g_2$ do not occur at the minima of $g_1$. Under this assumption we see that $\max \{ g_2 (t) \} < \max \{ g_1(t)\}$ if and only if each $t_{1, i , \max} \in (t_{g,i}, t_{g_, i+1})$, where $t_{g,i}$ and $t_{g_,i+1}$ are two subsequent intersection points of the reponses $g_1$ and $g_2$ and $t_{1,i,\max}$ is the $i^{\text{th}}$ maximum of $g_1$. If this condition were not satisfied there would exist $t$ such that $g_2(t) < \lambda$ or $g_2 (t) > \gamma_1$. Thus we see to ensure that $g_2 < g_1$ at $n$ maxima followed by $n$ minima $g_1$ and $g_2$ must have $2n$ intersection points. If the maximum is not followed by a minimum then it requires one intersection point. By assumption, only the final maximum may not be followed by a minimum.

		 Suppose the first minimum occurs before the first maximum. Either $t_{g,0} < t_{1,1, \min}$ or $t_{g,0} \in (t_{1,1,\min}, t_{1,1,\max})$ i.e. $u_1 (t) = u_2(t)$ for all $t < t_{g,0} - \varepsilon$, $\varepsilon \in (0,t_{g,0})$. This reduces to the case where the first maximum of $g_1$ occurs before the first minimum of $g_1$with distinct inputs $v_1(t,N-1)$ and $v_2(t,N-1)$, of the form \eqref{eq:genu}.

		 Suppose $g_1$ has $N+1$ minima, and therefore $N$ maxima. By the above $t_{g,0} < t_{1,1,\min}$ thus from Lemma \ref{lem:extend} $2N -1$ intersection points remain but $2N$ intersection points are required as each maximum is followed by a minimum.

		 Conversely, suppose $g_1$ has $N+1$ maxima, and therefore $N$ minima, this implies $g_2$ must intersect $g_1$ at $2N+1$ points. But by Lemma \ref{lem:extend} there exist at most $2N$ intersection points. 
		 
		 Finally, suppose that $m$ of the $t_{g,i} = t_{1,i,\min}$. This would imply that at least one of the maxima require one fewer intersection points. Indeed, $n$ maxima followed by $n$ minima would require $2n-m$ intersection points. At such $t_{g,i}$ we have that $h_1 = h_2$ and $g_1 = g_2$. Each intersection point after which $g_2 < g_1$ corresponds to a pulse of $u_2(t,N)$. Thus each $t_{1,i,\max}$ corresponds to a pulse of $u_2(t,N)$. Either there are $N+1$ maxima of $g_1$, in which case $u_2 (t,N)$ does not have sufficiently many pulse inputs, or there are $N$ maxima and $N+1$ minima of $g_1$. In this case $g_1$ must start with a minimum. As above we may assume that $t_{g,0} < t_{1,1,\min}$, which reduces to the case of $N$ maxima followed by $N$ minima and only $N-1$ pulses of $u_2(t,N)$ which we may apply.
		 
		 Hence, there exists no such response $g_2$ and therefore no input $u_2(t,N)$ exists. Thus $u_1 (t)$ is optimal.
	 \qed \end{pf}
		 
	 \begin{rem}
		 The converse of Theorem \ref{thm:mult} does not hold for all $w$. As rapid changes in $w$ will outpace the response time of $h$. However the results may be applied over specified bounded intervals, each of which have different maxima, to find the optimal input for each interval.
	 \end{rem}

	 \begin{lem}
		 \label{lem:thelastone}
		 Suppose $a < b < \infty$. Then there exists $M \geq 1$ such that the minimised response $\max\{g(h_n, w)\} =\max\{g(h_{n+1}, w)\}$ for all $t \in [a,b]$ and for all $n \geq M$, where $h_n$ is the response to the input $u(t , n)$.
	 \end{lem}

 \section{Example}

\begin{figure}[H]
	 \begin{flushleft}
		 \includegraphics[width=10cm, height=9.7cm, keepaspectratio]{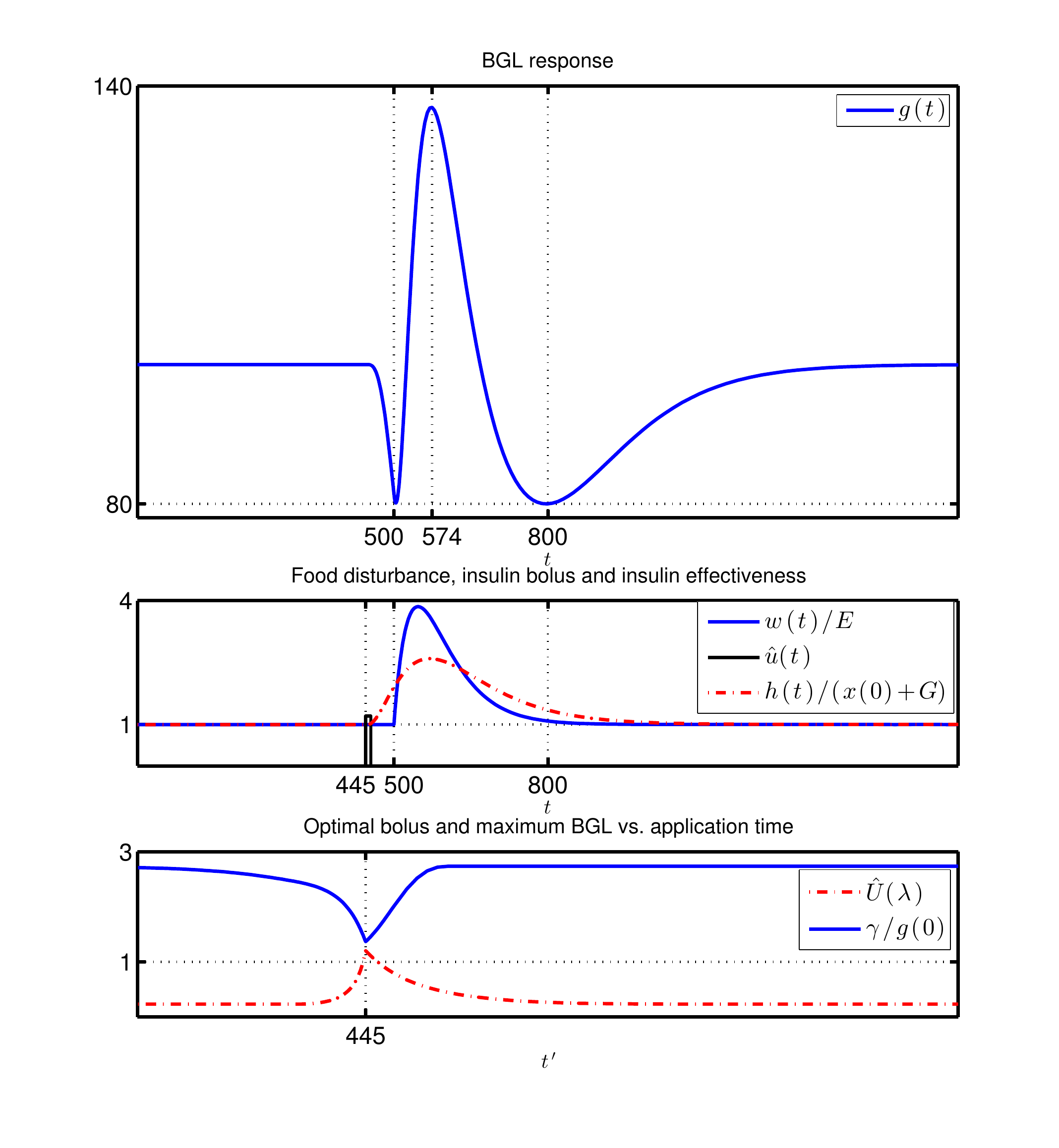}
	\end{flushleft}
	 \caption{The optimal glucose response to the functions $w(t)$ and $h(t)$, which are shown in the second plot, and the magnitude of a proper bolus and the maximum glucose concentration as a function of the delivery time $t'$, shown in the third plot.}
	 \label{fig:minmaxmin}
 \end{figure}

 In the example presented in Figure \ref{fig:minmaxmin}, we consider the following values for the parameters in \eqref{eq:eqs} and \eqref{eq:wandh}:
 $d=0.0204$, $k=497.5124$, $c=0.0213$, $a=0.0106$, $b=8.11 \times 10^{-4}$, $G=0.0032$, $E=1.3$, and $r(t) = 0.0018 f_1(t)$, where $f_1(t)$ is the response of the system of linear differential equations:
\[
	\begin{pmatrix}
		\dot f_1(t) \\
		\dot f_2(t)
	\end{pmatrix} 
	= 
	\l ( \frac{1}{47} \r)
	\begin{pmatrix}
		-1 &  1 \\
		0  & -1
	\end{pmatrix}
	\begin{pmatrix}
		f_1(t) \\
		f_2(t)
	\end{pmatrix} 
	+
	\begin{pmatrix}
		0 \\
		\delta(t)
	\end{pmatrix}
\]
where $\delta(t)$ is an impulse of magnitude $120$ applied at time $500$. We take the initial conditions to be as in Section \ref{sec:prelim} and set $g(\infty) = g(0) = 100 \mathrm{mgdl}^{-1} \, (5.2 \mathrm{mmolL}^{-1})$. The input is of the form \eqref{eq:u} with $\bar{u}$ computed as in \eqref{eq:basal} and duration $\tau=10$. The bolus magnitude $\hat u$ is computed as in Theorem~\ref{thm:bolus}, for the optimal injection time $t'$ that minimises $\gamma$, the global maximum of $g$. The minimum glucose concentration $\lambda$ is chosen to be $80 \mathrm{mgdl}^{-1} \, (4.4 \mathrm{mmolL}^{-1})$. 

 The first plot in Figure \ref{fig:minmaxmin} shows the plasma glucose response to the function, $w(t)$, shown in the second plot of the same Figure, and an optimal bolus input $\hat u$ delivered at time $445$. Two minima occur at times $500$ and $800$ bounding the unique maximum which occurs at time $574$. The final plot of Figure \ref{fig:minmaxmin} shows the maximum glucose concentration, $\gamma$, and the magnitude of a proper input bolus as a function of the input time $t'$. We see that $\gamma$ is minimised at the optimal input time $445$.

\section{Conclusions}

Current research aims to generalise the presented results to any bounded input function $u(t)$. We are also interested in studing the effect of $\tau$ -- the length of the bolus delivery interval $A$, on the response $g(t)$.

We do not know whether similar results to those presented may be shown for other models of glucose metabolism. In particular, those which include other factors such as exercise \citep{roy07} or free fatty acid metabolism \citep{roy06}. Given the general nature of the proofs of the current results we believe it is likely that similar results may hold for other models.

Additionally, we aim to derive a formula for the maximum plasma glucose concentration which is independent of the times $t_{\min}$ and $t_{\max}$ and depends solely on $\lambda$ and the set $A$ i.e. to find $f: \mathbb{R}_+ \times A \to \mathbb{R}_+$ such that $\gamma = f(\lambda,A)$. This may allow us to extend the results of Corollary \ref{cor:seesaw} by specifying the rate at which the maximum concentration increases with respect to increases in the fixed minimum concentration $\lambda$.

Finally, we desire to prove that there is an optimal partition of $\mathbb{R}_+$ into intervals so that the converse of Theorem~\ref{thm:mult} holds over each interval and no other partition will produce a lower maximum glucose concentration. Such a result may follow by extending Lemma \ref{lem:anotherone} to cover inputs of the form \eqref{eq:genu}. This would also allow us, in conjuction with Lemma \ref{lem:thelastone} to specify the minimum number of pulses required to achieve the lowest possible maximum glucose concentration over some bounded interval.

\bibliography{refs}
\end{document}